\documentclass[letterpaper, 10 pt, conference]{ieeeconf}  
\IEEEoverridecommandlockouts                              
\usepackage{amsfonts}
\usepackage{amssymb,amsmath}
\usepackage{graphicx}
\usepackage{cite}
\renewcommand{\epsilon}{\varepsilon}
\renewcommand{\phi}{\varphi}
\renewcommand{\theta}{\vartheta}

\newcommand{\diam}{\mathop{\rm diam}}
\renewcommand{\d}{\,\mathrm{d}}

\newcommand{\real}{\mathbb{R}}

\newtheorem{proposition}{Proposition}
\newtheorem{theorem}{Theorem}

\newtheorem{example}{Example}

\title{\LARGE \bf
Feedback maximum principle \\ for ensemble control of local continuity equations. \\ An application to supervised machine learning$^*$
}

\author{ Maxim Staritsyn$^{1}$, Nikolay Pogodaev$^{2}$, Roman Chertovskih$^{1}$, and
Fernando Lobo Pereira$^{1}$
\thanks{*The work is supported in part
by the Fundamental Research Program of the Russian State Academies of Sciences, Project
I.1.4.1; Russian Scientific Foundation, Projects 17-11-01093 and 19-11-00258; and FCT
(Portugal): R$\&$D Unit SYSTEC -- POCI-01-0145-FEDER-006933/SYSTEC funded by ERDF $|$
COMPETE2020 $|$ FCT/MEC $|$ PT2020 extension to 2018,  UID/\-EEA/\-00147/\-2020 and
NORTE-01-0145-FEDER-000033.}
\thanks{$^{1}$Maxim Staritsyn, Roman Chertovskih and Fernando Lobo Pereira are with 
Faculdade de Engenharia, Universidade do Porto, Rua Dr. Roberto Frias, s/n 4200-465, Porto, Portugal
        {\tt\small starmaxmath@gmail.com, roman@fe.up.pt, flp@fe.up.pt}
}%
\thanks{$^{2}$Nikolay Pogodaev is with Matrosov Institute for System Dynamics and Control Theory of SB~RAS, 134, Lermontov st., 664033, Irkutsk, Russia
        {\tt\small nickpogo@gmail.com}}%
}

\begin{document}
\maketitle
\thispagestyle{empty}
\pagestyle{empty}
\begin{abstract}
We consider an optimal control problem for a system of local continuity equations on a space of probability measures. Such systems can be viewed as macroscopic models of ensembles of non-interacting particles or homotypic individuals, representing several different ``populations''. For the stated problem, we propose a necessary optimality condition, which involves feedback controls inherent to the extremal structure, designed via the standard Pontryagin's Maximum Principle conditions. This optimality condition admits a realization as an iterative algorithm for optimal control. As a motivating case, we discuss an application of the derived optimality condition and the consequent numeric method to a problem of supervised machine learning via dynamic systems.
\end{abstract}
\section{INTRODUCTION}

This article is devoted to a particular (relatively simple) class of optimal control problems for transport equations in the space of probability measures. In the last few years, such systems
have become a very popular area of the modern applied mathematics, mainly thanks to the recent progress in the analysis on metric spaces of probability measures, and in the optimal transportation  theory \cite{AGS,SantambrogioBOOK,Villani}. Among prominent applications of such developments is the analysis of various  models of collective behavior, pedestrian dynamics and crowd motion control (see, e.g., the survey \cite{Piccoli} and citations therein). 

Our study is to some extent aligned with articles \cite{Chaudhari2018,Sonoda2019,WeinanJiequn}, which present another attractive field of application of differential equations in spaces of measures, namely, in the area of artificial neural networks (NNs); 
it was whipped up by recent works \cite{Benning,Cuchiero,Seidman,Weinan} that promote an approach to supervised machine learning via dynamical systems. This approach delegates the role of the ``machine'', supposed to solve a problem of classification or regression, to a controlled ODE, which can be interpreted as an ``infinitely deep'' NN.

To give the intuition, we treat the binary classification problem, where a machine is trained to ``weed out the bad from the good'', based on given sample data. The idea is, roughly speaking, to regard the set of samples 
as a crowd of individuals representing different ``populations'', and try to find a common control that steers these two populations to different targets.   
Imagine, for example, a hillock of wheat on a plate, 
contaminated by chaff. Can we separate the wheat from the chaff by shaking the plate in a horizontal plane? Theoretically, the answer is positive, if we idealize the cereals and the particles of chaff up to  material points $x_i \in \mathbb R^2$ and $y_j \in \mathbb R^2$, assume that they do not influence one another, and rightly choose the way of shaking.

Thanks to the known results on the controllability of ODE's implied by the classical Chow–Rashevskii theorem from the sub-Riemannian geometry, for any $T>0$, and any finite sets of points  $x_i \in \mathbb R^n$ and $y_j \in \mathbb R^n$, there is a finite number $m$ of vector fields $f_k$,\footnote{In fact, for any $n$, there are $m=5$ (!) desired vector fields  \cite{Cuchiero}.} and a smooth control $u=(u^1, \ldots, u^m)$ such that the respective solution $t \mapsto X_t[u](\eta)$  to the Cauchy problem 
\begin{equation*}
\label{int-ODE}\dot x = \sum_{k=1}^m f^k(x) \, u^k, \quad x(0)=\eta,
\end{equation*}
enjoys the property:
\[
X_T[u](x_i) = \xi_i, \quad X_T[u](y_j) = \zeta_j, \quad \forall i,j.
\]
Here, $\{\xi_i\} \subset \mathbb R^n$ and $\{\zeta_j\}\subset \mathbb R^n$ are finite sets of arbitrary pairwise different points having the cardinalities $N$ and $M$ 
of $\{x_i\}$ and $\{y_j\}$, respectively. 

In the above speculative example,  $\{\xi_i\}$ and $\{\zeta_j\}$ can be grouped arbitrarily close to two distant points, $\xi \in \mathbb R^2$ and $\zeta \in \mathbb R^2$, which means that the cereals and the chaff are indeed separated into two hillocks, i.e., they are \emph{classified}.

In practice, to provide
the above \emph{generic  universal  interpolation} property, one can take seven vector fields $f^k$, which mimic 
residual NNs \cite{Cuchiero}.  As soon as one has a control at hand with such a property, 
it can be applied to new data, beyond the training sample, and, if the points of these new data were 
close enough to some used in the 
training, it is natural to expect an accurate classification due to the continuity of the map $\eta\mapsto X_{(\cdot)}[u](\eta)$.

The beautiful result on the controllability, however, does not hint, how to design the actual control, which 
solves our problem: 
\[
\sum_{i} \left(X_T[u](x_i)-\xi\right)^2 + \sum_{j} \left(X_T[u](y_j)-\zeta\right)^2 \to \min.
\]
If the numbers of points $x_i$ and $y_j$ were few, we could recruit standard tools of the optimal control theory (Pontryagin's Maximum Principle and Dynamic Programming), as well as existing numeric algorithms, see, e.g., \cite{Benning,Qianxiao}. However, if the cardinalities $N$ and $M$ of the sets $\{x_i\}$ and $\{y_j\}$ 
are relatively large, the classical methods are either not applicable anymore or numerically inefficient.

In this article, we propose a method 
for numeric implementation of such problems for the case 
$N,M \to \infty$. Our approach combines the ``averaging principle'' 
from the control of
multi-agent 
systems with the technique of ``feedback control variations'' involving the constructions of 
Pontryagin's Maximum Principle, adapted from 
\cite{Dykhta2014,Dykhta2014-1}. 

By the first principle, the sets $\{x_i\}_{i=1}^N$ and $\{y_j\}_{j=1}^M$ are replaced by empirical probability measures $\frac{1}{N} \sum_{i=1}^N \delta_{x_i}$ and $\frac{1}{M} \sum_{j=1}^M \delta_{y_j}$, while the respective system of $N+M$ control ODEs turns into a system of two PDEs, namely, local transport equations. For the resulted distributed control problem, we design a deterministic optimization technique, which uses the 
classical Pontryagin's Maximum Principle, albeit in a non-standard way, to generate 
feedback controls. As we show, this approach 
serves to escape local extrema, and therefore has a potential to find a global solution.

\section{PROBLEM STATEMENT}

\subsection{Notations}

Given a metric space $X$, we use the standard notation $C([0,T]; X)$ for the spaces of continuous maps $[0,T]\mapsto X$ (endowing it with the usual $\sup$-norm), and $L_p([0,T];\mathbb R^n)$, $p=1,\infty$, for the Lebesgue quotient spaces of summable and bounded measurable functions $[0,T]\mapsto \mathbb R^n$, respectively.

By $\mathcal P =\mathcal P(\mathbb R^n)$, we denote the set of probability measures on the vector space $\mathbb R^n$, and by $\mathcal P_1=\mathcal P_1(\mathbb R^n)$ the subset of $\mu \in \mathcal P$ having finite first moment, i.e., such that
\(
\displaystyle\mathrm m_1(\mu)\doteq
\int_{\mathbb R^n} |x| \, d\mu(x)<\infty.\label{moment}
\)
Recall that $\mathcal P_1$ is a  complete separable metric space when it is endowed with the 1-Kantorovich (Wasserstein) distance $$
W_1(\mu, \nu)\doteq \sup\left\{\int_{\mathbb R^n}\varphi \, d(\nu\!-\!\mu)\left| \ \begin{array}{c}
\varphi \in C(\mathbb R^n; \mathbb R), \\ {\rm Lip}(\varphi)\leq 1
\end{array}\right.\right\}.
$$
Here, ${\rm Lip}(\varphi)$ is the minimal Lipschitz constant of $\varphi$.  

Given $\mu \in \mathcal P$  and a Borel measurable map $F: \, \mathbb R^n \mapsto \mathbb R^n$, $F_\sharp \mu$ denotes the push-forward $\mu\circ F^{-1} 
$ of $\mu$ through $F$. 

$\mathcal L^n$ stands for the usual Lebesgue measure on $\mathbb R^n$. 

\subsection{Finite-Dimensional Model}

To simplify the presentation, we will consider the case of two classes of individuals, ``$x$'' and ``$y$'', whose states at time $t$ are represented by vectors $x_i(t), y_j(t) \in \mathbb R^n$, $i=\overline{1,N}$, $j=\overline{1,M}$ (hereinafter $i=\overline{1,N}$ abbreviates $i\in \{1, \ldots, N\}$). An essential hypothesis, we are to accept here, is that all $x_i$ are homotypic and so are all $y_j$. In other words, $x_i$ and $x_k$ are indistinguishable to us, though we keep to distinguish $x_i$ from $y_j$, for any $i$ and $j$. 

The representatives of both classes find themselves in a common environment, but may have different properties, which means that the action of the media on each class may be different. This action results in a drift of points $x_i$ and $y_j$ in accordance with different vector fields $f$ and $g$.

Given a fixed time interval $[0,T]$, the individuals $x_i$ and $y_j$ are aimed at minimizing the values $\ell_1(x_i(T))$ and $\ell_2(y_j(T))$, where the cost functions $\ell_{1,2}: \, \mathbb R^n \to \mathbb  R$ are the same for all representatives of a class.

Playing the part of a guide, we can influence the media (i.e., \emph{all} the individuals of both classes, simultaneously) via a vector field $v: \, \mathbb R^n\times U \to  \mathbb R^n$
depending on a control parameter $u$ 
from a given set $U \subset \mathbb R^m$. As feasible control signals, we admit (equivalence classes of) 
measurable functions $[0,T] \mapsto U$. The space $\mathcal U \doteq L_\infty([0,T]; U)$ of 
controls is 
endowed with the weak-* topology $\sigma(L_\infty, L_1)$.

Thus, we come to the optimal control problem $(P_{N,M})$:
\begin{align}
&\mbox{inf}\left(\frac{1}{N}\sum_{i=1}^N\ell_1\left(x_i(T)\right) + \frac{1}{M}\sum_{j=1}^M\ell_2\left(y_j(T)\right)\right) 
\mbox{ subject to }  \nonumber\\
&  x_i(0)=x_i^0, \quad   \dot x_i=f(x_i)+ v_u(x_i), \quad i=\overline{1,N},
\label{ODE-x}\\
&   y_j(0)=y_j^0, \quad  \dot y_j=g(y_j)+ v_u(y_j), \quad j=\overline{1,M},
\label{ODE-y}\\
&t \in [0,T], \quad u \in \mathcal U,\label{u}
\end{align}
where $x_i^0$ and $y_j^0$ are given initial positions of the individuals. Note that $(P_{M,N})$ can be viewed as a version of the mean field social control problem (see, e.g., \cite{Huang,Wang}).

\subsection{Standing Assumptions}

We make the 
standard regularity hypotheses $(H_1)$: 

    $U \subset \mathbb R^m$ is compact; 
    vector fields $F_u(x) \doteq  f(x) + v_u(x)$ and $G_u(x) \doteq g(x) + v_u(x)$ are continuous in $(x,u)$, and there exist $C>0$ such that \[
    |F_u(x)| \leq C(1 +|x|)\mbox{ and }|F_u(x)-F_u(y)| \leq C |x-y|
    \]
    for any $u \in U$ and all $x, y \in \mathbb R^n$ (similarly for $G_u$);
    $\ell_{1,2}$ are continuous.

\subsection{Problem Reformulation in Terms of Probability Measures}

Now we shall pass to the so-called mean field limit of the ODE \eqref{ODE-x} and \eqref{ODE-y}. For this, consider the curves \(t \mapsto \mu^N_t \doteq \frac{1}{N} \sum_{i=1}^N \delta_{x_i(t)}\) and \(t \mapsto \nu^M_t\doteq \frac{1}{M} \sum_{j=1}^M \delta_{y_j(t)}\) in $\mathcal P_1(\mathbb R^n)$. Assume that
\begin{equation}
\mu^N_0 = \frac{1}{N} \sum_{i=1}^N \delta_{x_i^0}\rightharpoondown \vartheta_1 \mbox{ and }\nu^M_0 = \frac{1}{M} \sum_{i=1}^M \delta_{y_j^0}\rightharpoondown \vartheta_2\label{initMu}
\end{equation}
for some $\vartheta_{i} \in \mathcal P_1$, $i=1,2$.\footnote{$\rightharpoondown$ indicates the weak convergence of measures.} Then the curves $t \mapsto \mu^N_t$ and $t \mapsto \nu^M_t$ converge\footnote{This convergence follows from \eqref{initMu} and the continuous dependence of the solutions to \eqref{PDE-mu} and \eqref{PDE-nu} on the initial measures $\vartheta_{i}$, $i=1,2$ (see, e.g., \cite[Lemma~2.8]{Pogo-ContEq}).} in \(C([0,T]; \mathcal P_1)\) 
to distributional solutions $t \mapsto \mu_t$ and $t \mapsto \nu_t$ of the continuity equations
\begin{gather}
\partial_t \mu_t + \nabla \cdot\left(\left[f+v_{u(t)}\right] \, \mu_t\right)=0, \quad  \mu_0=\vartheta_1,\label{PDE-mu}\\
\partial_t \nu_t + \nabla \cdot\left(\left[g+v_{u(t)}\right] \, \nu_t\right)=0, \quad \nu_0=\vartheta_2.\label{PDE-nu}
\end{gather}
Here, $\partial_t$ abbreviates $\frac{\partial}{\partial t}$, 
``$\cdot$'' denotes the scalar product, and $\nabla \doteq \frac{\partial}{\partial x}$. Note that equations \eqref{PDE-mu} and \eqref{PDE-nu} are independent and paired only by the control function $u$. As it is standard, solutions to PDEs \eqref{PDE-mu} and \eqref{PDE-nu} are understood in the sense of distributions, and can be represented as
\[
\mu_t = (X_t)_\sharp \vartheta_1, \quad\nu_t = (Y_t)_\sharp \vartheta_2, \quad t \in [0,T],
\]
where $t \mapsto X_t(x)\doteq X_t[u](x) \in C(\mathbb R^n; \mathbb R^n)$ and $t \mapsto Y_t(x)$ denote the flows of the characteristic ODEs \eqref{ODE-x} and \eqref{ODE-y}.

Thus, we come to the optimal control problem $(P)$ in the space of probability measures:
\[
\inf_{u \in \mathcal U}\left(\int\ell_1 \d\mu_T[u]+\int\ell_2 \d\nu_T\left[u\right]\right) \mbox{ subject to \eqref{PDE-mu} and \eqref{PDE-nu}.}
\]
Hereinafter, we abbreviate $\int = \int_{\real^d}$ and drop arguments of integrands for brevity.

\subsection{Pontryagin's Maximum Principle}

To approach problem $(P)$, we will recruit the following assertion, which is a version of the paradigmatic 
Pontryagin's Maximum Principle (PMP).
\begin{proposition}\label{PMP}
Assume that hypotheses $(H_{1})$ hold together with the following additional assumption
\begin{description}
    \item[$(H_2)$:] 
    $F_u$ and $G_u$ 
    are  continuously differentiable in $x$. 
\end{description}
Let $\sigma\doteq (\mu, \nu, u)$ be an \emph{optimal control process}. Then, for $\mathcal L^1$-almost all (a.a.) $t \in [0,T]$, $u$ satisfies the \emph{maximum condition}
\begin{equation}\label{extCond}
u(t) \in
\arg\max_{u \in U} \left(\int \nabla p \cdot v_u \, \d \mu_t + \int \nabla q \cdot v_u \, \d \nu_t\right),
\end{equation}
where
$p$ and $q$ are solutions of the {dual transport equations}
\begin{gather}
\label{p}
\partial_t p + \nabla p \cdot (f + v_{u(t)})=0, \quad p_T=-\ell_1,\\
\label{q}
\partial_t q + \nabla q \cdot (g + v_{u(t)})=0, \quad q_T=-\ell_2.
\end{gather}
\end{proposition}
\smallskip

Solutions to \eqref{p}, \eqref{q} are also understood in the sense of distributions and are constructed by the method of characteristics as
\(
p_t = -\ell_1\circ(X_t)^{-1},\) \(q_t = -\ell_2\circ (Y_t)^{-1},\) \(t \in [0,T].\)

The proof of Proposition~\ref{PMP} follows closely the one in \cite[Theorem 2]{PogoStar2019}. This assertion
gives a necessary condition for \emph{local} extremum within the class of so-called needle variations of control function. Here, one can point out that $(P)$ is a \emph{nonlinear} (and nonconvex) problem, even if $v_u$ is linear in $u$: although the cost is linear
in the measure, the dynamics contains the product $u \mu$. This means that PMP shall not be a sufficient optimality condition. In other words, as it is inherited from control of ODEs, Proposition~\ref{PMP} does not select out \emph{global} solutions to problem $(P)$, and, therefore, it has a potential for ``improvement''. Towards this, we propose to adopt 
the approach \cite{Dykhta2014,Dykhta2014-1} from the finite-dimensional control theory, which enables to extract  additional information 
from the extremal condition \eqref{extCond}. 

The idea is to use the ``local information'' (a reference control process), to generate feedback signals of the ``extremal structure''. Realization of such feedbacks through a standard sampling scheme leads to 
a class of control variations, which 
appears to be somewhat richer than the usual class of needle variations. The optimality qualification within this new class of variations results in a 
necessary condition 
that remains within the formalism of PMP, but demonstrates a greater potential to discard non-optimal 
extrema.

\section{FEEDBACK MAXIMUM PRINCIPLE}
To simplify the exposition of the main idea, we now focus on the case of a single population by setting $\ell_1=\ell$, $\ell_2 \equiv 0$, $\vartheta_1=\vartheta$. 
Then, our problem reduces to
\[
I[u] \doteq \int \ell \d \mu_T \to \min_{u \in \mathcal U} \mbox{ subject to \eqref{PDE-mu}}.
\]

Similarly to \cite{PogoStar2019}, one can derive the following exact formula for the increment of the cost functional in terms of solutions to \eqref{PDE-mu} and \eqref{p}:
\begin{align*}
\Delta_u I =&  \int_0^T\d t\int \nabla \bar p_t \cdot \left(v_{\bar u(t)}-v_u\right) \, \d \mu_t.
\end{align*}
This formula entails the specification of the control function $u$ in the ensemble-feedback form
\begin{align}\label{Feedback}
u_t[\mu] \in U_t^{\bar p}(\mu)\doteq \arg\max_{u \in U}  \int \nabla \bar p_t(x) \cdot v_u(x) \, \d \mu(x),
\end{align}
in which case one could expect an ``improvement'' of the reference control $\bar u$:
\(
\Delta_u I < 0.
\)
However, 
the realization of \eqref{Feedback} as a control law requires the same accuracy as in the case of ODEs, due to the, possibly, discontinuous character of the ensuing vector fields.%

\subsection{Ensemble-Feedback Controls and Sampling Solutions}

In this section, we shall discuss the concept of solution to feedback controlled 
PDE \eqref{PDE-mu}. 
Consider two types of mappings, \(\mathbf{u}: \, [0,T]\times\mathcal P \to U\) and \(\mathfrak{u}: \, [0,T]\times\mathcal P \to \mathcal U,\)
that could be used to construct sampling solutions of the respective transport equations (to be discussed below); $\mathbf u$ could be arbitrary, while $\mathfrak u$ is assumed to be Borel measurable. We denote the sets of such functions by $\mathbf U$ and $\mathfrak U$, respectively.  Feedbacks of the class $\mathbf U$ generate piecewise constant controls in accordance with the classical Krasovskii-Subbotin scheme, while the ones of the class $\mathfrak U$ produce piecewise open-loop controls (short-term programs) in the spirit of the so-called model predictive control.

Now, we are going to present two sampling schemes, entailed by the proposed notions of feedback control, which are designed by usual Euler polygons.

\subsubsection{$\mathbf u$-Sampling Scheme}
Given $\mathbf u \in \mathbf U$, and 
a partition \(\pi=\{t_k\}_{k=0}^K\)  of the 
interval $[0,T]$,
\(
    t_0=0, \ t_{k-1}<t_k, \ k=\overline{1, K}, \ t_K=T,
\)
 the 
 polygonal arc \(t \mapsto \mu_t^\pi[\mathbf u]\in \mathcal P\) is defined via step-by-step integration of PDE 
 \eqref{PDE-mu} over 
 $[t_{k-1}, t_k]$: 
\begin{align}\label{sample-0}
    \mu^\pi_{t}[\mathbf{u}]=\mu^k_t, \quad t \in [t_{k-1}, t_{k}), \quad k =\overline{1,K},
 \mbox{ where}\\
   \label{sample-1}
    \mu^0\equiv \vartheta, \quad \mu^k\doteq \mu\big[\mathbf{u}_{t_{k-1}}[\mu^{k-1}_{t_k}]\big](t_{k-1},\mu^{k-1}_{t_k}),
    \end{align}
    and $\mu[u](\tau,\vartheta)$ is 
    a distributional solution of the continuity equation with control $u$, and the initial condition
    \(\mu_{\tau}=\vartheta\).

    Along with the polygonal arc, we compute a piecewise constant control
    \begin{equation}\label{sampleControl-1}
    u^\pi \doteq  \mathbf{u}_{t_{k-1}}[\mu^{k-1}] \mbox{ on } [t_{k-1}, t_{k}), \ k =\overline{1,K}.
    \end{equation}

Now, $\mathbf u$-sampling solution is introduced as 
any \emph{partial}  limit in $C([0,T]; \mathcal P_1)$ 
of a sequence of the above 
polygons
 as ${\rm diam}(\pi)\doteq \displaystyle\max_{1\leq k\leq K}(t_{k}-t_{k-1}) \to 0$; $\mathfrak S_1(\mathbf{u})$ stands for the the set of all such solutions, produced by $\mathbf u$.

 \smallskip

 \subsubsection{$\mathfrak u$-Sampling Scheme}
For $\mathfrak u \in \mathfrak U$, we design a polygonal arc \(t \mapsto \mu_t^\pi[\mathfrak u]\in \mathcal P\) by step-by-step integration \eqref{sample-0}, where $\mu^k$, $k=\overline{1,K}$, are 
defined by
   \begin{equation}\label{sample-2}
    \mu^0\equiv \vartheta, \quad \mu^k\doteq \mu\big[\mathfrak{u}_{(\cdot)}[\mu^{k-1}_{t_k}]\big](t_{k-1},\mu^{k-1}_{t_k}),
    \end{equation}
    and the associated control is a 
    concatenation 
    \begin{equation}\label{sampleControl-2}
    u^\pi(t) \doteq  \mathfrak{u}_{t}[\mu^{k-1}], \ t \in [t_{k-1}, t_{k}), \ k =\overline{1,K}.
    \end{equation}
Up to the mentioned difference, the notion of $\mathfrak u$-sampling solution and the set $\mathfrak S_2(\mathfrak{u})$ are introduced as above.

\smallskip

The following assertion is implied by a generalized version of the Arzela-Ascoli theorem, if we notice that, thanks to $(H_1)$, polygonal arcs $\mu_t^\pi\in C([0,T];\mathcal P_1)$ are uniformly Lipschitz continuous with a constant depending only on $C$, $T$ and $\mathrm m_1(\vartheta)$ \cite[Lemma 3]{Pogodaev2019} (and therefore, by the Gronwall's lemma, are equicontinuous, and uniformly bounded).

\smallskip

\begin{proposition}
Under assumptions $(H_1)$, $\mathfrak S_1(\mathbf{u})\neq \emptyset$ and $\mathfrak S_2(\mathfrak{u})\neq \emptyset$, for any $\mathbf{u} \in \mathbf U$ and $\mathfrak{u} \in \mathfrak U$.
\end{proposition}

\smallskip
Simple examples \cite[Appendix A]{StarSor2020} show that $\mathfrak S_1(\mathbf{u})$ and $\mathfrak S_2(\mathfrak{u})$ are not proper subsets of each other, even in the case of ODEs.

Given $\mathbf u$ or $\mathfrak u$, along with respective sampling solutions, one can consider the associated classical feedback solutions, defined analogously to the Carath\'{e}odory feedback solutions of ODEs, \cite{Dykhta2014}. Notice that, in general, there could be no classical feedback solutions at all.

\subsection{Optimality Qualification via Feedbacks
}
We are ready to present a necessary optimality condition, which employs feedback controls of the PMP-extremal structure, proposed by formula \eqref{Feedback}. For brevity, we will operate with feedbacks of the type $\mathbf U$ though it is possible to take the ones of the sort $\mathfrak U$, or both.

Consider a reference process $\bar \sigma = (\bar \mu, \bar u)$, whose optimality is to be checked. Taking $w \in \mathbf U$, let $\mathfrak S(w)$ denote the union of the set $\mathfrak S_1(w)$, and the, possibly empty, set of the respective classical feedback solutions. Introduce the following \emph{accessory} problem associated to $\bar \sigma$:
\[
(AP_{\bar\sigma})\ \, \mbox{Minimize}\!\int\ell \d \xi_T, \ \xi \in \mathfrak S(w), \ w=w_t[\mu] \in U^{\bar p}_t(\mu).
\]
Here, $U^{\bar p}_t$ is defined as in \eqref{Feedback}. Note that $(AP_{\bar\sigma})$ is a variational problem, and does not involve a control anymore. The following assertion presents a variational necessary optimality condition for a single-population version of problem $(P)$, which we call the \emph{feedback 
maximum principle} (FMP). 

\smallskip

\begin{theorem}\label{Th1}
  Assume that $(H_{1})$ and $(H_{2})$ hold. Then the optimality of $\bar \sigma = (\bar{\mu}, \bar u)$ for $(P)$ implies the optimality of $\bar \mu$  for $(AP_{\bar \sigma})$.
\end{theorem}
\smallskip
The proof 
is a combination of two simple facts:

i) an optimal process $\bar \sigma$ is PMP extremal, and therefore $\bar\mu$ is admissible for $(AP_{\bar \sigma})$ as a classical feedback solution (note that $\bar u(t) \in  U_t^{\bar p}(\bar \mu_t)$ for $\mathcal L^1$-a.a. $t \in [0,T]$), and


ii) sampling solutions can, by definition, be uniformly approximated by usual distributional solutions corresponding to open-loop controls.

Indeed, if one is managed to find a feedback control $\hat w \in U^{\bar p}$, and a sampling solution $\hat \xi \in \mathfrak  S_1(\hat w)$ such that
\(\int\ell \d \bar \mu_T> \int\ell \d \hat\xi_T,\)  
the same inequality also holds for polygonal arcs $\hat \xi^\pi$ approximating $\hat \xi$ and corresponding to certain piecewise constant open-loop control $\hat u^\pi \in \mathcal U$, for sufficiently small ${\rm diam}(\pi)$:
  \( \int\ell \d \bar \mu_T> \int\ell \d \hat \xi^\pi_T=\int \ell \d \mu[u^\pi],\)
which contradicts the optimality of $\bar\sigma$.

One has to admit that, in the presented variational form, the result is not easy to grasp. To benefit from FMP, one can use its counter-positive version, what means to treat Theorem~\ref{Th1} as a \emph{sufficient condition for non-optimality} via the qualification condition:
 \begin{equation}
  \int\ell \d \bar \mu_T \leq \int\ell \d \xi_T \quad \forall \xi \in \mathfrak S_1(w) \quad \forall w\in U^{\bar p}.
  \label{FMP}
  \end{equation}
In other words, given a reference process $\bar \sigma =(\bar \mu, \bar u)$ (which could be a local extremum) with a dual state $\bar p$, one can try to ``improve'' $\bar \sigma$ by testing it against feedbacks, running through the set $U^{\bar p}$.

In the following simple (but eloquent) example, PMP can not distinguish the ``worst'' process from the ``best'' one, while condition \eqref{FMP} does by ``playing'' with the non-uniqueness of signals of the extremal structure.

\smallskip

\begin{example}[Weeding out PMP extremum by FMP]
Consider a version of our problem in $\mathbb R^2$ with 
$\ell(a,b)=b$, 
\(
f(a,b)\equiv 0\), \(v_u(a,b)=(u, -a u),\)
and the following data:
\(
T =1\), \(U=[-1,1]\), and \(\vartheta(a,b) = \delta_{0}(a) \otimes \eta(b)\) ($\otimes$ denotes the product measure)
with an arbitrary $\eta\in \mathcal P(\mathbb R)$.

As one can easily check, the constant control $\bar u \equiv 0$ is PMP-extremal with the co-state
\(
\bar p_t(a,b) \equiv - b.
\)
The extremal multifunction is specified as
\(
U_t^{\bar p}(\mu) = {\rm sign} \, \int a \, \d \mu(a,b);
\)
${\rm sign} \, s= s/|s|$ for $s \neq 0$, and  ${\rm sign} \, 0= [-1,1]$ by definition.

Let us apply Theorem~\ref{Th1}. Since \(\int a  \d \vartheta=\int a  \d\delta_0(a)=0,\) we can choose any control $u \in {\rm sign} \, 0= [-1,1]$ as the initial guess. In our option, take $u = 1$. This control shifts the mass from the axis $a=0$ to the right ($a>0$); both sampling schemes leave the strategy $u=1$ intact and lead to a single solution
\(
t \mapsto \mu_t= \left[(a,b)\mapsto (a+t, b-at-t^2/2)\right]_\sharp \vartheta,
\)
which gives
\(
\int b \d \eta(b) - \frac{1}{2} = I[u] < I[\bar u] = \int b \d \eta(b).
\)

One can check that this is, in fact, a global solution, and the choice $u= -1$ leads to the same result. Furthermore, any $u\in [-1,1]\setminus\{0\}$ 
also improves the reference extremal $(\bar u, \bar \mu)$, and, therefore, qualifies it as non-optimal.
\end{example}

\subsection{Iterations of Feedback Maximum Principle.
Strategy for Numeric Implementation}

FMP can be implemented iteratively in the following manner, taken desired accuracies $\epsilon_{i}>0$, $i=1,2$:
\begin{enumerate}
    \item[0.] Given an initial guess $\bar u \in \mathcal U$, compute $\bar p$ as a solution to \eqref{p} with $u=\bar u$.

    \item[1.] Choose a feedback control
    $w \in U^{\bar p}$. 
    \item[2.] Fixed a partition $\pi \subset [0,T]$, compute the polygonal arc $\mu^\pi[w]$ 
    together with the open-loop 
    control
    \(
    u^\pi.
    \)
    \item[3.] Check the ``improvement'' property
    \(
    I[u^\pi]-I[\bar u]<\epsilon_1.
    \)
    If it holds then update the control:
    \(
    \bar u=u^\pi.
    \)
    Otherwise, check
    \(
    {\rm \diam}(\pi)<\epsilon_2.
    \)
    If ``yes'', return to step 2. Else, go to step 1.

\end{enumerate}

In practice, control variations provided by feedbacks $w \in U^{\bar p}$ could be too strong, such that the associated feedback solution would not ensure the improvement property. To fix this problem, on Step~2, one can use 
a ``localized'' feedback
\begin{equation}
\omega_t[\mu]=\alpha \, \bar u(t) + (1-\alpha) \, w_t[\mu]
\label{update}
\end{equation}
involving extra parameter $\alpha \in (0,1]$ to be adjusted.


Let us stress two features of the proposed conceptual numeric method. Notice that the control at each step is designed based on the aggregated, macroscopic information, which is an integral characteristic of the ensemble of 
particles. 
In practice, this means that 
there is no need to trace actual positions of all individuals, as it is enough to follow their ``averaged representative'', and this
provides an essential dimensionality reduction.

Another point to be mentioned is that a realization of the proposed algorithm requires numeric solution of Cauchy problems only. The latter is essentially simpler from the computational viewpoint than a shooting method for a boundary-value problem, or dynamic programming via 
the Hamilton-Jacobi equation, which is known to be a more complicated task than the original optimal control problem. 

\section{MACHINE LEARNING}\label{Sec:Ex}

Now, let us return to the multi-population framework and revise the correlation problem discussed in the Introduction.

For the binary case (an extension to several classes is trivial) this problem is exactly of the form $(P_{N,M})$, if we regard the collections of training samples  associated to different labels as representatives of different populations, specify the dynamics as
\(
f=g\equiv 0\), \(
v_u(x)=\sum 
f^k(x) \, u^k,\)
with $f^k: \, \mathbb R^n \to \mathbb R^n$ providing the generic  universal  interpolation, and take $\ell_1(x)=|x-\xi|^2$, $\ell_2(x)=|x-\zeta|^2$ with two (distant) target points $\xi\neq \zeta \in \mathbb R^n$. Thus, all the results of this paper are applicable to the classification problem.

Above, we regarded the set of training samples as a crowd of individuals. However, a sample itself could be viewed as a population. 
In the following example, we learn a machine to classify simple geometric forms 
having only \emph{single representative of each a class} as a training set. Such problems of  ``machine learning over small data'', which sound deliberately pointless within the 
paradigm of NNs, still can be dealt with by the discussed control-theoretical approach. For this, planar figures are treated as distributions on a plane, and the type of a figure marks out the respective population. 
If the data are properly prepared, the ``classification effect'' provided by a control is naturally extrapolated from a single sample measure (one copy of a figure) to all $W_1$-close measures (sufficiently similar figures) thanks to the continuity of the map $(\vartheta, u) \mapsto \mu[u](\vartheta)$. 
\begin{example}[Control system as a classifier for MNIST] Consider a very standard NN exercise: to design a classifier 
of handwritten digits, using the canonical database MNIST.
To our preference, 
monochromatic 
pictures representing 
``3'' and ``6'' were selected. As ``right answers'', two sample images (one for ``3'' and one for ``6'') were randomly chosen. The pictures are treated as discrete measures 
$\frac{1}{N} \sum_{i=1}^N\delta_{x_i}$ on the unit square, where $N$ and $x_i$ correlate with the power of white color. To learn the resulted finite-dimensional control system, a problem of the form $(P_{N,M})$ was solved numerically, taken
$m=3$, $\xi=-\zeta=(3,0)$, $f^1\equiv (0,1)$,
\(
f^2(a, b)=(a+b, a)\), \(f^3(a, b)=(\sin(a), a-b)\), and $U=\left\{(u^1, u^2, u^3): \, |u^{1,2,3}| \leq 1\right\}$. The obtained control (Fig.~\ref{control}) appeared to provide 
an accurate classification in a series of tests with randomly chosen candidates.

\begin{figure}[thpb]
      \centering
      \includegraphics[
      scale=0.4
      ]{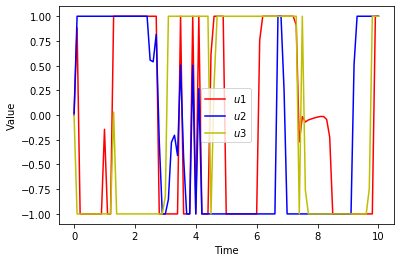}
      \caption{Numeric control for distinguishing ``3'' and ``6'' in MNIST.
      \label{control}
      }
  \end{figure}
\end{example}

Above, the problem was still finite-dimensional. In the next example, our algorithm separates two \emph{infinite} populations which constitute geometric patterns (a ``cross'' and a ``ring'') corrupted by a Gaussian blur. 

\begin{example}[Classification of ``uncertain'' planar figures] 

The same $f^k$ and $U$ are used as in Example~2. 
Taken $\xi=(1,1)$, $\zeta=-\xi$, and $T=2$, Cauchy problems 
(\ref{PDE-mu}), (\ref{PDE-nu}) and (\ref{p}), (\ref{q}) were solved in the spatial computational domain $(a, b)\in[-5,5]^2$ with Dirichlet boundary conditions corresponding to the vanishing measures on the boundaries. 
A standard upwind scheme (the donor-cell upwind method \cite{leveque}) is used, backward for (\ref{p}), (\ref{q}) and forward in time for (\ref{PDE-mu}), (\ref{PDE-nu}). The mesh spacing in both spatial directions is 
$h=5\cdot 10^{-2}$, the time step is 
$\tau=2\cdot 10^{-3}$, $\alpha$ in (\ref{update}) is chosen to be $0.75$.   

The initial measures $\mu_0$ and $\nu_0$ are depicted in Fig.~\ref{munuini}, representing a cross and a ring, correspondingly. 
The initial guess for the control is constant: $\bar u 
\equiv (1, 0, 0)$. After 30 iteration of the algorithm, 
we assume that the problem is solved. 
The terminal distributions $\mu_T$, $\nu_T$, and the respective control $u=(u^1, u^2,
u^3)$ are shown in Figs~\ref{munufin}, and~\ref{figu}.  

\begin{figure}
\vspace{0.4cm}
      \centering\includegraphics[scale=0.61]{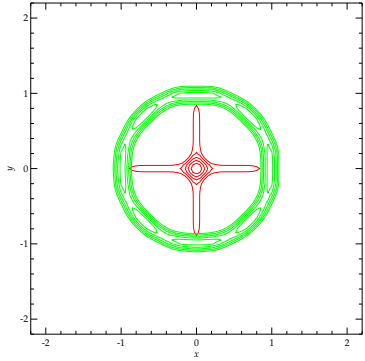}
      \caption{Isolines of $\mu_0$ (red) and $\nu_0$ (green). Isovalues are 0.9, 0.8, 0.7, 0.6 and 0.5 of the corresponding maximal value.\label{munuini}}
  \end{figure}
 
\begin{figure}
      \centering\includegraphics[scale=0.61]{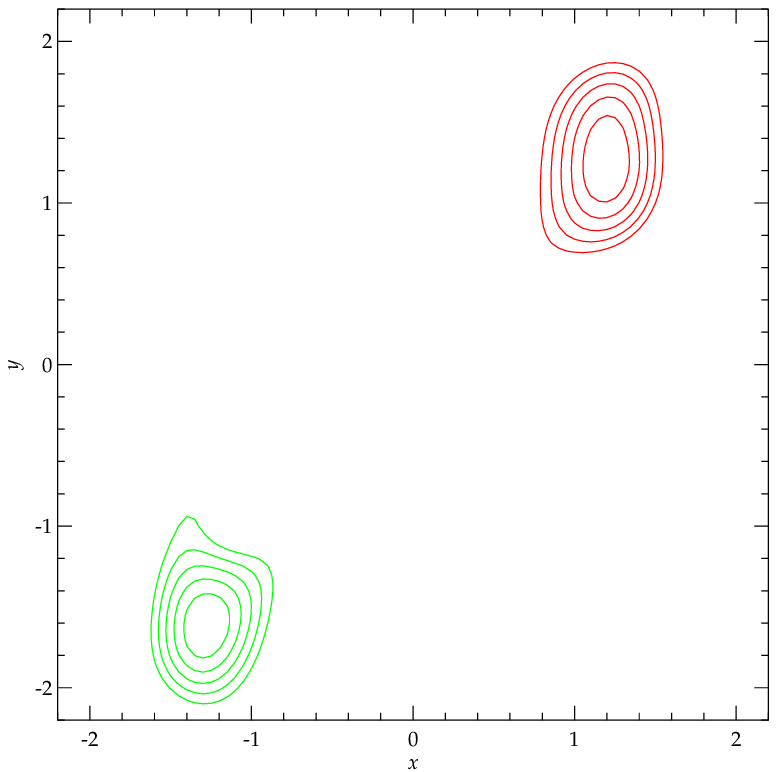}
      \caption{Isolines of $\mu_T$ (red) and $\nu_T$ (green). Isovalues are 0.9, 0.8, 0.7, 0.6 and 0.5 of the corresponding maximal value.\label{munufin}}
\end{figure}

\begin{figure}
      \centering\includegraphics[scale=0.87]{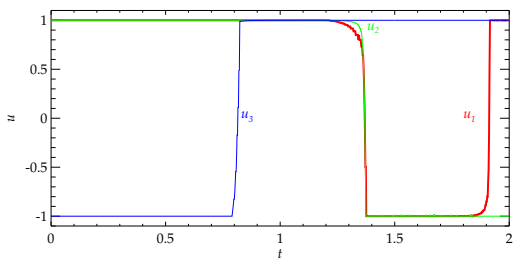}
      \caption{``Optimal'' control components: $u^1$ (red), $u^2$ (green) and $u^3$
      (blue).\label{figu}}
\end{figure}

\end{example}

\section{CONCLUSION}

A natural direction of future work 
would be an extension of the discussed framework to the case of nonlocal vector fields 
involving a convolutional term $K*\mu$, which brings certain analogy with convolutional NNs. The PMP for ensemble control of such nonlinear transport 
equations does exist \cite{POGODAEV20203585}, but the Hamiltonian system appears to be inseparable into the direct and dual subsystems, which makes the realization of our approach a challenging issue.

In what concerns applications to machine learning, the real potential of the presented approach needs investigation via numeric experiments, and this is another part of our future study. At present, as a 
mission of control-theoretical considerations of machine learning tasks, 
we see the new light shed
on 
mathematical principles behind the machinery of the artificial intelligence. 

\bibliographystyle{IEEEtran}
\bibliography{IEEEabrv,Staritsyn}

\end{document}